\documentclass{amsart}

\usepackage{epsfig}

\def\N{{\mathbb N}}
\def\C{{\mathbb C}}

\def\Z{{\mathbb Z}}
\def\d{{\rm d}}
\def\deg{{\rm deg\,}}
\def\deg{{\rm deg\,}}

\begin{document}

\title{\bf On polynomials orthogonal 
to all powers of a Chebyshev polynomial on a segment.
}
\author {F. Pakovich}

\maketitle

\markboth{}{}

\section{Introduction.}

In the recent series of papers \cite{bfy1}-\cite{bfy5} of M. Briskin,
J.-P. Francoise and Y. Yomdin the following ``polynomial moment problem"
was proposed as an infinitesimal version of the center problem
for the Abel differential equation in the complex domain:
{\it for a
complex polynomial $P(z)$ and distinct $a,b \in \C$
such that $P(a)=P(b)$ to describe polynomials $q(z)$
such that} $$\int^b_a P^i(z)q(z)\d z=0
\eqno(1)$$
{\it for all integer non-negative $i$.}

The following ``composition condition" imposed on $P(z)$ and $Q(z)=\int
q(z)\d z$ is sufficient
for polynomials $P(z),q(z)$ to satisfy (1):
{\it there exist polynomials $\tilde P(z),$
$\tilde Q(z),$ $W(z)$ such that} $$P(z)=\tilde P(W(z)),
\ \ Q(z)=\tilde Q(W(z)), \ {\it and} \
W(a)=W(b).\eqno(2)$$
Indeed, the sufficiency of condition (2) is a direct corollary of the
Cauchy theorem since after the change of variable
$z \rightarrow W(z)$ the new way of integration is closed.
It was suggested in the papers cited above (``the composition
conjecture") that condition (2) is actually equivalent to condition (1).
This conjecture was verified in several special cases. In particular,
when $a,b$ are not critical points of $P(z)$ (\cite{c}),
when $P(z)$ is indecomposable (\cite{pa2}),
and in some other special cases (\cite{bfy1}-\cite{bfy5},
\cite{ro}, \cite{pry}). Nevertheless, in general the composition
conjecture is not true.

A class of counterexamples to the composition conjecture was constructed in \cite{pa1}. The simplest of them has the following form:
$$P(z)=T_6(z), \ \ q(z)=T_3^{\prime}(z)+T_2^{\prime}(z),\ \ a=-\sqrt 3/2,\ \ b=\sqrt 3/2, $$ where $T_n(z)=\cos (n\arccos
z)$ is the n-th Chebyshev polynomial.
Indeed, since $T_2(\sqrt 3/2)=T_2(-\sqrt 3/2)$ it follows from the equality 
$T_6(z)=T_3(T_2(z))$ that (1) is satisfied for $P(z)=T_6(z)$ and $q_1(z)=T_2^{\prime}(z).$ 
Similarly, from $T_6(z)=T_2(T_3(z))$ and $T_3(\sqrt 3/2)=T_3(-\sqrt 3/2)$
one concludes that (1) holds for $P(z)=T_6(z)$ and $q_2(z)=T_3^{\prime}(z).$ Therefore, by linearity,
condition (1) is satisfied also for $P(z)=T_6(z)$ and $q(z)=q_1(z)+q_2(z).$
Nevertheless, for $P(z)=T_6(z)$ and $Q(z)=T_3(z)+T_2(z)$ 
condition (2) does not hold. 

More generally, it was shown in \cite{pa1} that 
any polynomial ``double decomposition" $A(B(z))=C(D(z))$  
such that $B(a)=B(b),$ $D(a)=D(b)$
supplies counterexamples to the composition conjecture whenever $\deg B(z),$ $\deg D(z)$ are coprime. 
Note that double decompositions with $\deg A(z)=\deg D(z),$ $\deg B(z)=\deg C(z)$ and $\deg B(z),\deg D(z)$
coprime are described explicitly by Ritt's theory of factorization
of polynomials. They are equivalent either to decompositions with $A(z)=z^nR^m(z),$
$B(z)=z^m,$ $C(z)=z^m,$ $D(z)= z^nR(z^m)$ for a polynomial $R(z)$
and $(n,m)=1$ or to decompositions with $A(z)=T_m(z),$ $B(z)= T_n(z),$ $C(z)= T_n(z),$ $D(z)=
T_m(z)$
for Chebyshev polynomials $T_n(z),$ $T_m(z)$ and $(n,m)=1$ (see \cite{ri},
\cite{sh}).

In this paper we give a solution of the polynomial moment problem (1)
in the case when $P(z)$ is a Chebyshev polynomial $T_n(z).$
Denote by $V(T_n,a,b)$ the vector space 
over $\C$ consisting of polynomials $q(z)$ satisfying (1) for $P(z)=T_n(z).$ 
Note that any polynomial $T_m^{\prime}(z)$ such that 
$T_d(a)=T_d(b)$ for $d={\rm GCD}(n,m)$ is contained in $V(T_n,a,b)$ since
$T_{n}(z)=T_{n/d}(T_{d}(z))$ and $T_{m}(z)=T_{m/d}(T_{d}(z)).$ 

\vskip 0.2cm

\noindent{\bf Theorem.} {\it 
For any $n\in \N$ and $a,b\in \C$ such that $T_n(a)=T_n(b)$ polynomials
$T_m^{\prime}(z)$ such that $T_d(a)=T_d(b)$ for $d={\rm GCD}(n,m)$ form a
basis of $V(T_n,a,b)$. }

\vskip 0.2cm

The theorem implies that 
for Chebyshev polynomials the following weakened version of the composition conjecture is true: if $q(z)$ satisfies condition (1) with $P(z)=T_n(z)$ then 
$\int q(z)\d z$ can be represented as a {\it sum} of
polynomials for which condition (2) is hold. Moreover, actually the number of terms in 
such a representation can be reduced to two.

\vskip 0.2cm

\noindent{\bf Corollary.} {\it 
Let $q(z)\in V(T_n,a,b).$ Then  
there exist divisors $d_1,d_2$ of $n$ such that 
$\int q(z) \d z =A(T_{d_1}(z))+B(T_{d_2}(z))$ for some polynomials $A(z),B(z)$
and
the equalities $T_{d_1}(a)=T_{d_1}(b),$ $T_{d_2}(a)=T_{d_2}(b)$
hold.

}

\vskip 0.2cm

For instance, if a polynomial $q(z)$ belongs to $V(T_6,-\sqrt{3}/2,\sqrt{3}/2)$ then 
the polynomial $\int q(z) \d z$ can be represented as $\int q(z) \d z =A(T_3(z))+B(T_2(z))$ 
for some polynomials $A(z),B(z).$
Note that such a representation in general is not unique
in contrast to the one for $q(z)$ providing by the theorem.

\section{Proofs.}
\subsection{Reduction.}
First of all, we establish that the theorem can be reduced to the
following
statement: {\it if $q(z)=Q^{\prime}(z)$ is contained
in $V(T_n,a,b)$ then} 
$$T_{d}(a)=T_{d}(b) \ \ \ {\it for} \ \ \ d={\rm GCD}(n,\deg Q).\eqno(3)$$
Indeed, assuming that this statement is true the theorem can be deduced 
as follows. For $q(z)\in V(T_n,a,b)$
set $m_0=\deg Q(z)$ and define $c_{0}\in \C$ by the condition that the degree of $Q_{1}(z)=Q(z)-c_{0}T_{m_0}(z)$ is strictly less then 
$m_0.$ Since for $d_0={\rm GCD}(n,m_0)$ the equalities
$$T_{n}(z)=T_{n/d_0}(T_{d_0}(z)), \ \ \ \ T_{m_0}(z)=T_{m_0/d_0}(T_{d_0}(z))$$
hold it follows from $T_{d_0}(a)=T_{d_0}(b)$ that $T_{m_0}^{\prime}(z)\in V(T_n,a,b).$
Therefore, by linearity, $Q_1^{\prime}(z)\in V(T_n,a,b).$
If $\deg Q_1(z)=m_1$ then, similarly, for some $c_{m_1}\in \C$ we have
$Q_1(z)=c_{m_1}T_{m_1}(z)+Q_2(z),$ where 
$Q_2^{\prime}(z)\in V(T_n,a,b)$ and $\deg Q_2(z) < m_1.$ 

Continuing in the same way  
and observing that $m_{i+1}<m_{i}$
we eventually arrive to the representation 
$$\int q(z)\d z = \sum_{i=0}^{k}c_iT_{m_i}(z), \ \ \ c_i\in \C,\eqno(4)$$ such that $T_{d_i}(a)=T_{d_i}(b)$ for
$d_i={\rm GCD}(n,m_i).$ Since polynomials of different degrees are
linearly independent over $\C$ we conclude that the polynomials
$T_m^{\prime}(z)$ such that $T_d(a)=T_d(b)$ 
for $d={\rm GCD}(n,m)$ form a
basis of the vector space $V(T_n,a,b).$

\subsection{Proof of the theorem for non-singular $a,b$.}
By 2.1 it is enough to show that condition 
(1) with $P(z)=T_n(z),$ $q(z)=Q^{\prime}(z)$ implies condition (3). On the
other hand, it is known (see \cite{c} or
\cite{pa2}) that for any polynomial $P(z)$ such that $a,b$ are not
critical points of
$P(z)$ the conditions
(1) and (2) are equivalent. Therefore, it is enough to prove that
(2) with $P(z)=T_n(z)$ implies (3). 

Suppose that (2) holds and set $w=\deg W(z).$ 
Since by Engstrom's theorem (see e.g. \cite{sh}, Th. 5)
for any double decomposition $A(B(z))=C(D(z))$ we have
$$[\C(B,D):\C(D)]=\deg D/{\rm GCD}(\deg B,\deg D),$$ it follows from the
equality 
$$T_n(z)=\tilde P(W(z))=T_{n/w}(T_{w}(z))$$ 
that $\C(W)=\C(T_w)$. Therefore, since $W(z), T_w(z)$ are polynomials, there exists a linear
function $\sigma(z)$ such that $W(z)=\sigma(T_w(z))$ and, hence,
$W(a)=W(b)$ yields 
$T_{w}(a)=T_{w}(b).$ Since $w$ is a divisor of 
$d={\rm GCD}(n,\deg Q)$ the decomposition 
$T_{d}(z)=T_{d/w}(T_{w}(z))$ holds and, therefore, $T_{w}(a)=T_{w}(b)$
implies
$T_{d}(a)=T_{d}(b).$

\subsection{Necessary condition for $P(z), q(z)$ to satisfy (1).}
To investigate the case when at least one from the points $a,b$ is 
a critical point of $T_n(z)$ we will 
use a condition, obtained 
in \cite{pa2}, formula (6) and in a more general situation in \cite{pry}, 
Theorem 3.9, which
is necessary for polynomials $P(z),$ $q(z)$ to satisfy (1).
To formulate this condition let us introduce the following notation.
Say that a domain $U\subset \C$
is admissible with respect to the polynomial $P(z)$ if $U$ is
simply connected and contains no
critical values of $P(z).$ By the monodromy
theorem, in such a domain there exist $n=\deg P(z)$ single-valued
branches of $P^{-1}(z).$ Let $U$  
be an admissible with respect to $P(z)$ domain
such that its boundary $\partial U$ contains the point $z_0=P(a)=P(b).$
Denote by $p_{u_1}^{-1}(z),$
$p_{u_2}^{-1}(z), ... , p_{u_{d_a}}^{-1}(z)$
(resp. $p_{v_1}^{-1}(z),$ $p_{v_2}^{-1}(z), ... ,
p_{v_{d_b}}^{-1}(z)$)
the branches of $P^{-1}(z)$ in $U$ which map points close to $z_0$
to points close to the point $a$ (resp. $b$). In particular, the number
$d_a$ (resp. $d_b$) equals the multiplicity of the point $a$ (resp. $b$)
with respect to $P(z).$ 
In this notation 
a necessary condition for $P(z), q(z)$ to satisfy (1) 
has the following form: 
{\it if $P(z),$ $q(z)=Q^{\prime}(z)$ satisfy (1) then 
in any admissible with respect to $P(z)$ domain $U$ such that 
$z_0\in \partial U$ the equality}
$$d_b\sum_{s=1}^{d_a} Q(p_{u_s}^{-1}(z)) \equiv d_a\sum_{s=1}^{d_b}
Q(p_{v_s}^{-1}(z)) \eqno(5)$$ {\it holds.}

\subsection{Monodromy of $T_n(z)$.}
To make condition (5) useful we must examine
the monodromy group 
of $T_{n}(z).$ It follows from 
$T_n(\cos\phi)=\cos(n\phi),$ $n\geq 1,$ that
finite critical values of polynomial $T_n(z)$ are $\pm 1$
and that the preimages of $\pm 1$ are
$\cos (\pi j /n),$ $j=0,1, ... , n.$
To visualize the monodromy group of $T_{n}(z)$
consider the preimage $\lambda=P^{-1}[-1,1]$ of the segment $[-1,1]$ under the map $P(z)\,:\,\C\rightarrow \C.$
It is convenient to consider $\lambda$ as a bicolored graph
$\Omega$ embedded into the Riemann sphere.
By definition, white (resp. black) vertices of $\lambda$ are preimages of the point $1$  
(resp. $-1$) and edges of $\lambda$ are preimages of the interval $(-1,1)$.
Since the multiplicity of each critical point of $T_n(z)$ equals 2,
the graph $\lambda$ is a ``chain-tree''
and, as a point set in $\C,$ coincides with the segment $[-1,1]$
(see fig. 1).
\begin{figure}[h] 
\centering
\input{fche2.pstex_t}
\caption{ }
\end{figure}

Let us fix an admissible with respect to $T_n(z)$ domain $U$ such
that $U$ is unbounded and 
contains the interval $(-1,1).$ Any branch 
$T_{n,j}^{-1}(z),$ $0\leq j \leq
n-1,$ of $T^{-1}_n(z)$ in $U$ maps the interval  
$(-1,1)$ onto an edge of $\lambda$ and we will label such an edge by
the symbol $l_j$ (an explicit numeration of the branches of $T^{-1}_n(z)$
will be defined later).
Denote by $\pi_{1}\in S_n$
(resp. $\pi_{-1}, \pi_{\infty}\in S_n$) the
permutation defined by the condition that
the analytic continuation of the
functional element $\{U,T_{n,j}^{-1}(z)\},$ $0\leq j \leq n-1,$
along a clockwise
oriented loop around $1$ (resp. $-1,$ $\infty$) is 
the functional element $\{U,T_{n,\pi_{1}(j)}^{-1}(z)\}$
(resp. $\{U,T_{n,\pi_{-1}(j)}^{-1}(z)\},$
$\{U,T_{n,\pi_{\infty}(j)}^{-1}(z)\}$).
The tree $\lambda$ represents the monodromy group of $T_{n}^{-1}(z)$
in the following sense:
the edges of $\lambda$ are identified with branches of $T_{n}^{-1}(z)$
and the permutation $\pi_{1}$ (resp. $\pi_{-1}$) is 
identified with the permutation arising under 
clockwise rotation of edges of $\lambda$ around white (resp. black) vertices.\footnote{Note that any polynomial with two finite critical values can be represented by an appropriate bicolored plane tree
and vice versa; it is a very particular case of the Grothendieck 
correspondence between Belyi functions and graphs embedded
into compact Riemann surfaces (see e.g. \cite{Schn}).}

In order to fix a convenient numeration of branches of $T_n^{-1}(z)$ in $U$ consider an auxiliary domain $U_{\infty}=U\cap B,$
where $B$ is a disc with the center at the infinity 
such that branches of $T_n^{-1}(z)$ can be represented in $B$
by their Puiseux expansions at infinity.
In more details, if $z^{\frac{1}{n}}$ denotes a fixed branch of the
algebraic
function which is inverse to $z^n$ in $U_{\infty},$ then each branch of
$T^{-1}_n(z)$ can be represented in $U_{\infty}$  
by the convergent series 
$$\phi_j(z)=\sum_{k=-\infty}^{1}t_k\varepsilon_n^{jk}z^{\frac{k}{n}}, \ \ \ 
t_k\in \C, \ \ \ \varepsilon_n=exp(2\pi i/n), \eqno(6)$$
for certain $j,$ $0\leq j \leq n-1.$
Now we fix a numeration of branches of $T_n^{-1}(z)$ in $U$ as follows: the branch $T_{n,j}^{-1}(z),$ $0\leq j \leq n-1,$ is the analytic continuation of 
$\phi_j(z)$ 
from $U_{\infty}\cap U$ to $U$ and the branch $z^{\frac{1}{n}}$
is defined by the condition that
$T_{n,0}^{-1}(z)$ maps the interval  
$(-1,1)$ onto the interval $(\cos (\pi /n),1).$
Since the result of the analytic continuation of the
functional element $\{U_{\infty}, \varepsilon_n^{j}z^{\frac{1}{n}}\},$ $0\leq j \leq n-1,$ along a clockwise
oriented loop around $\infty$ is the functional element $\{U_{\infty}, \varepsilon_n^{j+1}z^{\frac{1}{n}}\},$ such a choice of the numeration implies that $\pi_{\infty}=(0\, 1\, 2\, ...\, n-1).$
Furthermore, it follows from
$\pi_{\infty}\pi_{-1}\pi_{1}=1,$ 
taking into account the combinatorics of $\lambda,$
that the numeration of edges of $\lambda$ coincides with the one 
indicated on the figure 1 that is $\pi_{-1}=(0\,\, n\!-\!1)(1\, \, 
n\!-\!2)(2\, \,
n\!-\!3)
... $ 
and $\pi_{1}=(1\, \, n\!-\!1)(2\, \, n\!-\!2)(3\, \, n\!-\!3) ... $.

\subsection{Proof of the theorem for singular $a,b$.} 
Again, it is enough to establish that (3) holds.  
Let $Q^{\prime}(z)\in V(T_n,a,b)$ with $\deg Q(z)=m.$
Since at least one from points $a,b$ is a critical point
of $T_n(z)$ the number $z_0=P(a)=P(b)$ equals $\pm 1.$
Suppose at first that $z_0=1.$ Then
$a=\cos (2j_1\pi /n),$ $b=\cos ( 2j_2\pi /n)$ 
for certain 
$j_1,j_2,$ $0\leq j_1,j_2 \leq [n/2],$ 
and 
condition (5) has the following form:
$$Q(T_{n,j_1}^{-1}(z))+Q(T_{n,n-j_1}^{-1}(z))
=Q(T_{n,j_2}^{-1}(z))+Q(T_{n,n-j_2}^{-1}(z)), \eqno(7)$$
where $T_{n,i}^{-1}(z)$ is represented in 
$U_{\infty}$ by the series (6). 
Since $t_{1}\neq 0,$ the 
comparison of the leading coefficients 
of the corresponding (7) Puiseux expansions gives
$$\varepsilon_n^{j_1m}+\varepsilon_n^{(n-j_1)m}
=\varepsilon_n^{j_2m}+\varepsilon_n^{(n-j_2)m}.$$ Therefore, the
number $\varepsilon_n^{\frac{m}{d}},$ where $d={\rm GCD}(n,m),$ is a root
of the
polynomial
$$f(z)=z^{j_1d}+z^{(n-j_1)d}-z^{j_2d}-z^{(n-j_2)d}.$$
Since $\varepsilon_n^{\frac{m}{d}}$ is a primitive $n$-th root of unity 
and the coefficients of $f(z)$ are integers, this fact implies that
$n$-th cyclotomic polynomial $\Phi_n(z)$ divides $f(z)$ in the ring
$\Z[z]$ and, therefore, that
the primitive $n$-th root of unity 
$\varepsilon_n$ also
is a root of $f(z)$. 
Hence, 
$$\varepsilon_n^{j_1d}+\varepsilon_n^{-j_1d}
=\varepsilon_n^{j_2d}+\varepsilon_n^{-j_2d}.$$   
Since 
$$a=\cos (2j_1\pi
/n)=\frac{1}{2}(\varepsilon_n^{j_1}+\varepsilon_n^{-j_1}),
\ \ \ b=\cos
(2j_2\pi /n)=\frac{1}{2}(\varepsilon_n^{j_2}+\varepsilon_n^{-j_2}),$$
it follows now from $$T_d(\frac{1}{2}(z+\frac{1}{z}))
=\frac{1}{2}(z^d+\frac{1}{z^d}) \eqno(8)$$ that $T_{d}(a)=T_{d}(b).$

Similarly, if $z_0=-1,$ 
assuming that $a=\cos ((2j_1+1)\pi /n),$ $b=\cos ((2j_2+1)\pi /n)$
for certain $j_1,j_2,$ $0\leq j_1,j_2 \leq [(n-1)/2],$ 
we obtain the equality 
$$T_{n,j_1}(z)+T_{n,n-j_1-1}(z)=T_{n,j_2}(z)+T_{n,n-j_2-1}(z)$$
which implies 
$$\varepsilon_{n}^{j_1m}+\varepsilon_{n}^{(n-j_1-1)
m}
=\varepsilon_{2n}^{j_2m}+\varepsilon_n^{(n-j_2-1)
m}$$ and 
$$\varepsilon_n^{j_1d}+\varepsilon_n^{-(j_1+1)d}
=\varepsilon_n^{j_2d}+\varepsilon_n^{-(j_2+1)d}.$$
It yields that
$$\varepsilon_{2n}^{2j_1d}+\varepsilon_{2n}^{-2(j_1+1)d}
=\varepsilon_{2n}^{2j_2d}+\varepsilon_{2n}^{-2(j_2+1)
d},$$ where $\varepsilon_{2n}=exp (2\pi i/2n),$
and, multiplying the last eqaulity by $\varepsilon_{2n}^d,$ we get
$$\varepsilon_{2n}^{(2j_1+1)d}+\varepsilon_{2n}^{-(2j_1+1)d}
=\varepsilon_{2n}^{(2j_2+1)d}+\varepsilon_{2n}^{-(2j_2+1)
d}.$$ Since
$$a=\frac{1}{2}(\varepsilon_{2n}^{2j_1+1}+\varepsilon_{2n}^{-(2j_1+1)}), \ \ \ b=\frac{1}{2}(\varepsilon_{2n}^{2j_2+1}+\varepsilon_{2n}^{-(2j_2+1)}),$$
we conclude as above that $T_{d}(a)=T_{d}(b).$

\subsection{Proof of the corollary.} Suppose $q(z)\in V(T_n,a,b).$ Then,
by the theorem, $\int q(z) \d z $ can be represented by the sum (4).  We will prove
the corollary by induction on the number of non-zero terms in this
representation. Since for each $i,$ $0\leq i \leq k,$ in (4) we have
$T_{m_i}(z)=T_{m_i/d_i}(T_{d_i}(z))$ with $d_i={\rm GCD}(n,m_i)$ and 
$T_{d_i}(a)=T_{d_i}(b),$ the corollary
is true for $k=0,1.$ 

Suppose now that $k>1.$ By the inductive assumption
there exist $r,s$ such that $$\sum_{i=0}^{k-1} c_iT_{m_i}(z)=
A(T_r(z))+B(T_s(z))$$ and $T_r(a)=T_r(b),$ $T_s(a)=T_s(b).$ Choose $v,w\in
\C$ such that $$a=\frac{1}{2}(v+\frac{1}{v}), \ \ \
b=\frac{1}{2}(w+\frac{1}{w}).$$ Then, by (8),
$T_r(a)=T_r(b)$ implies that
$$\frac{1}{2}(v^{r}+\frac{1}{v^{r}})=\frac{1}{2}(w^{r}+\frac{1}{w^{r}}).$$
In its turn, the last equality holds if and only if $v^{r}=w^{\mu_1 r},$
where $\mu_1=\pm 1.$ Similarly, the equalities $T_s(a)=T_s(b),$
$T_{m_k}(a)=T_{m_k}(b)$ yield that $v^{s}=w^{\mu_2 s},$ $v^{m_k}=w^{\mu_3
m_k},$ where $\mu_2,\mu_3=\pm 1.$

Suppose $\mu_1=\mu_2=\mu.$ Then
$v^{r}=w^{\mu r},$ $v^{s}=w^{\mu s}$ 
and an easy reasoning involving 
roots of unity shows that
$d={\rm GCD}(r,s)>1$ and $v^{d}=w^{\mu d}.$ Therefore, by
(8), $T_d(a)=T_d(b)$ 
and we can represent $\int q(z)\d z$ as $$\int q(z)\d z
=C(T_d(z))+c_kT_{m_k}(z),$$ where 
$C(z)=A(T_{r/d}(z))+B(T_{s/d}(z)).$ 

If $\mu_1=-\mu_2$ then either 
$\mu_1=\mu_3$ or $\mu_2=\mu_3$ and we conclude as above that either 
$$\int q(z)\d z =E(T_e(z))+B(T_s(z)),$$ where 
$E(z)=A(T_{r/e}(z))+c_kT_{m_k/e}(z),$ $e={\rm GCD}(r,m_k),$ and
$T_e(a)=T_e(b)$ or
$$\int q(z)\d z =A(T_r(z))+F(T_f(z)),$$ where 
$F(z)=B(T_{s/f}(z))+c_kT_{m_k/f}(z),$
$f={\rm GCD}(s,m_k),$ and $T_f(a)=T_f(b).$

\bibliographystyle{amsplain}

\begin{thebibliography}{10}



\bibitem {bfy1} M. Briskin, J.-P. Francoise, Y. Yomdin,
\textit{Une approche au probleme du centre-foyer de Poincare,}
C. R. Acad. Sci., Paris, Ser. I, Math. 326, No.11, 1295-1298 (1998).


\bibitem {bfy2} M. Briskin, J.-P.
Francoise, Y. Yomdin,
\textit{Center conditions, compositions of polynomials and moments on
algebraic curve}, Ergodic Theory Dyn. Syst. \textbf{19}, no 5,
1201-1220 (1999).


\bibitem {bfy3} M. Briskin, J.-P. Francoise, Y. Yomdin,
\textit{Center condition II: Parametric and model center problems}, Isr.J.
Math. \textbf{118}, 61-82 (2000).


\bibitem {bfy4} M. Briskin, J.-P. Francoise, Y. Yomdin,
\textit{Center condition III: Parametric and model center problems}, Isr.
J. Math. \textbf{118}, 83-108 (2000).

\bibitem {bfy5} M. Briskin, J.-P. Francoise, Y. Yomdin,
\textit{Generalized moments, center-focus conditions and compositions of
polynomials,} in "Operator theory, system theory and related topics",
Oper. Theory Adv. Appl.,  \textbf{123}, 161--185 (2001).

\bibitem {c} C. Christopher, \textit{Abel equations: composition
conjectures and
the model problem}, Bull. Lond. Math. Soc. \textbf{32}, No.3,
332-338 (2000).

\bibitem {pa1} F. Pakovich, \textit{A counterexample to the composition
conjecture}, Proc. Amer. Math. Soc. {\bf 130}, 3747-3749 (2002). 

\bibitem {pa2} F. Pakovich, \textit{On the polynomial moment
problem}, submitted.  


\bibitem {pry} F. Pakovich, N. Roytvarf and Y. Yomdin, \textit{Cauchy type
integrals of
Algebraic functions,} preprint.


\bibitem {ri} J. Ritt, \textit{Prime and composite polynomials,} Trans.
Amer. Math. Soc.  \textbf{23}, no. 1, 51--66 (1922)

\bibitem {ro} N. Roytvarf, \textit{Generalized moments, composition of
polynomials and Bernstein classes}, in "Entire functions in modern
analysis. B.Ya. Levin memorial volume", Isr. Math. Conf. Proc.
\textbf{15}, 339-355 (2001).


\bibitem {sh} A. Schinzel, \textit{Polynomials with special regard to
reducibility}, Encyclopedia of Mathematics and Its Applications
\textbf{77}, Cambridge University Press, 2000.


\bibitem {Schn} "The Grothendieck Theory of Dessins D'enfants"
(L. Shneps eds.),
Cambridge University Press ("London mathematical society
lecture notes series", vol. {\bf 200}), 1994.



\end{thebibliography}

\end{document}